\documentclass[reqno,a4paper]{amsart}
\usepackage{amssymb,isorot}
\usepackage[colorlinks=true,hypertex]{hyperref}
\date{}
\allowdisplaybreaks[4]
\theoremstyle{plain}
\newtheorem{thm}{Theorem}
\theoremstyle{remark}
\newtheorem{rem}{Remark}
\DeclareMathOperator{\td}{d\mspace{-2mu}}

\begin{document}

\title[The function $(b^x-a^x)/x$: Ratio's properties]
{The function $\boldsymbol{(b^x-a^x)/x}$: Ratio's properties}

\author[B.-N. Guo]{Bai-Ni Guo}
\address[B.-N. Guo]{School of Mathematics and Informatics,
Henan Polytechnic University, Jiaozuo City, Henan Province, 454010, China}
\email{\href{mailto: B.-N. Guo <bai.ni.guo@gmail.com>}{bai.ni.guo@gmail.com}, \href{mailto: B.-N. Guo <bai.ni.guo@hotmail.com>}{bai.ni.guo@hotmail.com}}

\author[F. Qi]{Feng Qi}
\address[F. Qi]{Research Institute of Mathematical Inequality Theory, Henan Polytechnic University, Jiaozuo City, Henan Province, 454010, China}
\email{\href{mailto: F. Qi <qifeng618@gmail.com>}{qifeng618@gmail.com}, \href{mailto: F. Qi <qifeng618@hotmail.com>}{qifeng618@hotmail.com}, \href{mailto: F. Qi <qifeng618@qq.com>}{qifeng618@qq.com}}
\urladdr{\url{http://qifeng618.spaces.live.com}}

\begin{abstract}
In the paper, after reviewing the history, background, origin, and applications of the functions $\frac{b^{t}-a^{t}}{t}$ and $\frac{e^{-\alpha t}-e^{-\beta t}}{1-e^{-t}}$, we establish sufficient and necessary conditions such that the special function
$\frac{e^{\alpha t}-e^{\beta t}}{e^{\lambda t}-e^{\mu t}}$
are monotonic, logarithmic convex, logarithmic concave, $3$-log-convex and $3$-log-concave on $\mathbb{R}$, where $\alpha,\beta,\lambda$ and $\mu$ are real numbers satisfying $(\alpha,\beta)\ne(\lambda,\mu)$, $(\alpha,\beta)\ne(\mu,\lambda)$, $\alpha\ne\beta$ and $\lambda\ne\mu$.
\end{abstract}



\thanks{The second author was partially supported by the China Scholarship Council}

\thanks{This paper was typeset using \AmS-\LaTeX}

\maketitle

\section{Introduction}

Recall~\cite{note-on-neuman.tex, note-on-neuman-ITSF-simplified.tex} that a $k$-times differentiable function $f(t)>0$ is said to be $k$-log-convex on an interval $I$ if
\begin{equation}\label{k-times-ineq}
0\le[\ln f(t)]^{(k)}<\infty,\quad k\in\mathbb{N}
\end{equation}
on $I$; If the inequality~\eqref{k-times-ineq} reverses then $f$ is said to be $k$-log-concave on $I$.
\par
For $b>a>0$, let
\begin{equation}\label{g0}
G_{a,b}(t)=
\begin{cases}
\dfrac{b^{t}-a^{t}}{t}, & t\neq 0; \\
\ln b-\ln a, & t=0.
\end{cases}
\end{equation}
In~\cite{jmaa-ii-97, (b^x-a^x)/x}, the complete monotonicity and inequality properties of $G_{a,b}(t)$ were first investigated. In~\cite{emv-log-convex-simple.tex, Cheung-Qi-Rev.tex, exp-funct-further.tex, exp-funct-appl-means.tex}, the $3$-log-convex and $3$-log-concave properties of $G_{a,b}(t)$ were shown. The function $G_{a,b}(t)$ has close relationships with the incomplete gamma function~\cite{schext, schext-rgmia, cubo}. It was ever used to prove the Schur-convex properties~\cite{schext, schext-rgmia, exp-funct-appl-means.tex}, the logarithmic convexities~\cite{Cheung-Qi-mean-rgmia, Cheung-Qi-mean, pams-62, pams-62-rgmia, Cheung-Qi-Rev.tex, exp-funct-appl-means.tex}, the monotonicity~\cite{ql, qx3} of the extended mean values. It was applied in~\cite{Gauchman-Steffensen-pairs, steffensen-qi-cheng-rgmia, steffensen-pair-Anal, qcw, onsp, onsp-rgmia} to construct Steffensen pairs. It was also employed in~\cite{notes-best.tex} to verify Elezovi\'c-Giordano-Pe\v{c}ari\'c's theorem~\cite[Theorem~1]{egp} which is related to the monotonicity of a function involving the ratio of two gamma functions. For more information, please refer to~\cite{bullen-handbook, cubo} and closely-related references therein.
\par
For $b>a>0$, let
\begin{equation}
F_{a,b}(t)=
\begin{cases}
\dfrac{t}{e^{bt}-e^{at}},& t\ne 0;\\[0.6em]
\dfrac1{b-a},& t=0.
\end{cases}
\end{equation}
In~\cite{Best-Constant-exponential.tex, best-constant-one-simple.tex, best-constant-one.tex, best-constant-rgmia, wuzh, wuzh-rgmia, best-constant-one-simple-real.tex}, \cite[p.~217]{souza} and~\cite[p.~295]{3rded}, the inequalities, monotonicity and logarithmic convexities of the function $F_{a,b}(t)$ for $a=b-1$ and its logarithmic derivatives of the first and second orders are established. In~\cite{best-constant-one.tex}, the history, background and origin of $F_{a,b}(t)$ for $a=b-1$ and its first two logarithmic derivatives were cultivated. In~\cite{remiander-Sen-Lin-Guo.tex, psi-reminders.tex, psi-reminders.tex-rgmia}, the logarithmic derivative of $F_{a,b}(t)$ for $a=b-1$ was applied to study the complete monotonicity of remainders of the first Binet formula and the psi function. In~\cite{ijmest-bernoulli, bernoulli-luo-guo-qi-rgmia, bernoulli-luo-guo-qi-debnath-IJMMS, euler-bernoulli-luo-qi-adv, bernoulli-qi-guo-rgmia}, the function $F_{\ln a,\ln b}(t)$ was utilized to generalize Bernoulli numbers and polynomials. In~\cite{emv-log-convex-simple.tex, exp-funct-further.tex}, the $3$-log-convex and $3$-log-concave properties of $F_{a,b}(t)$ were shown, among other things.
\par
For real numbers $\alpha$ and $\beta$ satisfying $\alpha\ne\beta$, $(\alpha,\beta)\ne(0,1)$ and $(\alpha,\beta)\ne(1,0)$, let
\begin{equation}
Q_{\alpha,\beta}(t)=
\begin{cases}
\dfrac{e^{-\alpha t}-e^{-\beta t}}{1-e^{-t}},&t\ne 0;\\
\beta-\alpha,&t=0.
\end{cases}
\end{equation}
In~\cite{mon-element-exp-final.tex, mon-element-exp.tex-rgmia, comp-mon-element-exp.tex, notes-best.tex}, the monotonicity and logarithmic convexities of $Q_{\alpha,\beta}(t)$ were discussed and the following conclusions were procured:
\begin{enumerate}
\item
The function $Q_{\alpha,\beta}(t)$ is increasing on $(0,\infty)$ if and only
if $(\beta-\alpha)(1-\alpha-\beta)\ge0$ and $(\beta-\alpha) (|\alpha-\beta|
-\alpha-\beta)\ge0$;
\item
The function $Q_{\alpha,\beta}(t)$ is decreasing on $(0,\infty)$ if and only
if $(\beta-\alpha)(1-\alpha-\beta)\le0$ and $(\beta-\alpha) (|\alpha-\beta|
-\alpha-\beta)\le0$;
\item
The function $Q_{\alpha,\beta}(t)$ is increasing on $(-\infty,0)$ if and only
if $(\beta-\alpha)(1-\alpha-\beta)\ge0$ and $(\beta-\alpha) (2-|\alpha-\beta|
-\alpha-\beta)\ge0$;
\item
The function $Q_{\alpha,\beta}(t)$ is decreasing on $(-\infty,0)$ if and only
if $(\beta-\alpha)(1-\alpha-\beta)\le0$ and $(\beta-\alpha) (2-|\alpha-\beta|
-\alpha-\beta)\le0$;
\item
The function $Q_{\alpha,\beta}(t)$ is increasing on $(-\infty,\infty)$ if and
only if $(\beta-\alpha) (|\alpha-\beta| -\alpha-\beta)\ge0$ and
$(\beta-\alpha) (2-|\alpha-\beta| -\alpha-\beta)\ge0$;
\item
The function $Q_{\alpha,\beta}(t)$ is decreasing on $(-\infty,\infty)$ if and
only if $(\beta-\alpha) (|\alpha-\beta| -\alpha-\beta)\le0$ and
$(\beta-\alpha) (2-|\alpha-\beta| -\alpha-\beta)\le0$;
\item
The function $Q_{\alpha,\beta}(t)$ on $(-\infty,\infty)$ is logarithmically
convex if $\beta-\alpha>1$ and logarithmically concave if $0<\beta-\alpha<1$;
\item
If $1>\beta-\alpha>0$, then $Q_{\alpha,\beta}(t)$ is $3$-log-convex on $(0,\infty)$ and $3$-log-concave on $(-\infty,0)$; if $\beta-\alpha>1$, then $Q_{\alpha,\beta}(t)$ is $3$-log-concave on $(0,\infty)$ and $3$-log-convex on $(-\infty,0)$.
\end{enumerate}
The monotonicity of $Q_{\alpha,\beta}(t)$ on $(0,\infty)$ was applied in~\cite{mon-element-exp-final.tex, sandor-gamma-3.tex, sandor-gamma-3-note.tex-aam, sandor-gamma-3-note.tex} to present necessary and sufficient conditions such that some functions involving ratios of the gamma and $q$-gamma functions are logarithmically completely monotonic. The logarithmic convexities of $Q_{\alpha,\beta}(t)$ on $(0,\infty)$ was used in~\cite{notes-best-new-proof.tex, notes-best.tex} to supply alternative proofs for Elezovi\'c-Giordano-Pe\v{c}ari\'c's theorem. For detailed information, please refer to~\cite{Wendel-Gautschi-type-ineq.tex, Wendel2Elezovic.tex} and related references therein.
\par
The functions $G_{a,b}(t)$, $F_{a,b}(t)$ and $Q_{\alpha,\beta}(t)$ have the following relations:
\begin{equation}
\begin{aligned}
G_{a,b}(t)&=\frac1{F_{\ln a,\ln b}(t)}, & F_{a,b}(t)&=\frac1{G_{e^b,e^a}(t)},\\
Q_{\alpha,\beta}(t)&=\frac{G_{e^{-\alpha},e^{-\beta}}(t)}{G_{1,e^{-1}}(t)},& Q_{\alpha,\beta}(t)&=\frac{F_{0,-1}(t)}{F_{-\alpha,-\beta}(t)}.
\end{aligned}
\end{equation}
\par
For real numbers $\alpha,\beta,\lambda$ and $\mu$ satisfying $(\alpha,\beta)\ne(\lambda,\mu)$, $(\alpha,\beta)\ne(\mu,\lambda)$, $\alpha\ne\beta$ and $\lambda\ne\mu$, let
\begin{equation}\label{ht-dfn}
H_{\alpha,\beta;\lambda,\mu}(t)=
\begin{cases}
\dfrac{e^{\alpha t}-e^{\beta t}}{e^{\lambda t}-e^{\mu t}},&t\ne 0,\\[0.5em]
\dfrac{\beta-\alpha}{\lambda-\mu},&t=0.
\end{cases}
\end{equation}
For positive numbers $r,s,u$ and $v$ satisfying $(r,s)\ne(u,v)$, $(r,s)\ne(v,u)$, $r\ne s$ and $u\ne v$, let
\begin{equation}\label{pt-dfn}
P_{r,s;u,v}(t)=
\begin{cases}
\dfrac{r^{t}-s^{t}}{u^{t}-v^{t}},&t\ne 0,\\[0.5em]
\dfrac{\ln r-\ln s}{\ln u-\ln v},&t=0.
\end{cases}
\end{equation}
It is clear that
\begin{equation}\label{ratio-relation-1}
H_{\alpha,\beta;\lambda,\mu}(t)=P_{e^\alpha,e^\beta;e^\lambda,e^\mu}(t)\quad\text{and}\quad  P_{r,s;u,v}(t)=H_{\ln r,\ln s;\ln u,\ln v}(t).
\end{equation}
In addition, the functions $H_{\alpha,\beta;\lambda,\mu}(t)$ and $P_{r,s;u,v}(t)$ can be represented as ratios of $G_{a,b}(t)$, $F_{a,b}(t)$ and $Q_{\alpha,\beta}(t)$ as follows:
\begin{equation}\label{ratio-relation}
H_{\alpha,\beta;\lambda,\mu}(t)=\frac{F_{\lambda,\mu}(t)}{F_{\alpha,\beta}(t)} =\frac{Q_{-\alpha,-\beta}}{Q_{-\lambda,-\mu}}\quad \text{and} \quad P_{r,s;u,v}(t)=\frac{G_{r,s}(t)}{G_{u,v}(t)}.
\end{equation}
\par
Since the functions $G_{a,b}(t)$, $F_{a,b}(t)$ and $Q_{\alpha,\beta}(t)$ have a long history, a deep background and many applications to several areas, so we continue to study the monotonicity and logarithmic convexities of their ratios, $H_{\alpha,\beta;\lambda,\mu}(t)$ and $P_{r,s;u,v}(t)$.
\par
Our main results may be stated as follows.

\begin{thm}\label{exp-gen-T1}
For real numbers $\alpha,\beta,\lambda$ and $\mu$ with $(\alpha,\beta)\ne(\lambda,\mu)$, $(\alpha,\beta)\ne(\mu,\lambda)$, $\alpha\ne\beta$ and $\lambda\ne\mu$, let
\begin{align}
\mathcal{A}&=(\alpha-\beta)(\alpha+\beta-\lambda-\mu), \\ \mathcal{B}&=(\alpha-\beta)(\alpha+\beta-|\alpha-\beta|-2\lambda),\\
\mathcal{C}&=(\alpha-\beta)(\alpha+\beta+|\alpha-\beta|-2\lambda),\\
\mathcal{D}&=(\alpha-\beta)(\alpha+\beta+|\alpha-\beta|-2\mu),\\
\mathcal{E}&=(\alpha-\beta)(\alpha+\beta-|\alpha-\beta|-2\mu).
\end{align}
Then the function $H_{\alpha,\beta;\lambda,\mu}(t)$ has the following properties:
\begin{enumerate}
\item
The function $H_{\alpha,\beta;\lambda,\mu}(t)$ is increasing on $(0,\infty)$ if and only if either $\lambda>\mu$, $\mathcal{A}\ge0$ and $\mathcal{C}\ge0$ or $\lambda<\mu$, $\mathcal{A}\le0$ and $\mathcal{B}\le0$.
\item
The function $H_{\alpha,\beta;\lambda,\mu}(t)$ is decreasing on $(0,\infty)$ if and only if either $\lambda<\mu$, $\mathcal{A}\ge0$ and $\mathcal{B}\ge0$ or $\lambda>\mu$, $\mathcal{A}\le0$ and $\mathcal{C}\le0$.
\item
The function $H_{\alpha,\beta;\lambda,\mu}(t)$ is increasing on $(-\infty,0)$ if and only if either $\lambda>\mu$, $\mathcal{A}\ge0$ and $\mathcal{E}\ge0$ or $\lambda<\mu$, $\mathcal{A}\le0$ and $\mathcal{D}\le0$.
\item
The function $H_{\alpha,\beta;\lambda,\mu}(t)$ is decreasing on $(-\infty,0)$ if and only if either $\lambda>\mu$, $\mathcal{A}\le0$ and $\mathcal{E}\le0$ or $\lambda<\mu$, $\mathcal{A}\ge0$ and $\mathcal{D}\ge0$.
\item
The function $H_{\alpha,\beta;\lambda,\mu}(t)$ is increasing on $(-\infty,\infty)$ if and only if either $\lambda>\mu$, $\mathcal{C}\ge0$ and $\mathcal{E}\ge0$ or $\lambda<\mu$, $\mathcal{B}\le0$ and $\mathcal{D}\le0$.
\item
The function $H_{\alpha,\beta;\lambda,\mu}(t)$ is decreasing on $(-\infty,\infty)$ if and only if either $\lambda>\mu$, $\mathcal{C}\le0$ and $\mathcal{E}\le0$ or $\lambda<\mu$, $\mathcal{B}\ge0$ and $\mathcal{D}\ge0$.
\item
The function $H_{\alpha,\beta;\lambda,\mu}(t)$ on $(-\infty,\infty)$ is logarithmically convex if $\frac{\alpha-\beta}{\lambda-\mu}>1$ or logarithmically concave if $0<\frac{\alpha-\beta}{\lambda-\mu}<1$.
\item
The function $H_{\alpha,\beta;\lambda,\mu}(t)$ is $3$-log-convex on $(0,\infty)$ and $3$-log-concave on $(-\infty,0)$ if either $\lambda-\mu>\alpha-\beta>0$ or $\alpha-\beta<\lambda-\mu<0$; The function $H_{\alpha,\beta;\lambda,\mu}(t)$ is $3$-log-concave on $(0,\infty)$ and $3$-log-convex on $(-\infty,0)$ if either $\alpha-\beta>\lambda-\mu>0$ or $\lambda-\mu<\alpha-\beta<0$.
\end{enumerate}
\end{thm}

\begin{thm}\label{exp-gen-T2}
For positive numbers $r,s,u$ and $v$ with $(r,s)\ne(u,v)$, $(r,s)\ne(v,u)$, $r\ne s$ and $u\ne v$, let
\begin{gather*}
\mathfrak{A}=\ln\frac{rs}{uv}\ln\frac{r}s,\quad \mathfrak{B}=\biggl(\ln\frac{rs}{u^2} -\biggl|\ln\frac{r}s\biggl|\biggr)\ln\frac{r}s, \quad \mathfrak{C}=\biggl(\ln\frac{rs}{u^2} +\biggl|\ln\frac{r}s\biggl|\biggr)\ln\frac{r}s,\\
\mathfrak{D}=\biggl(\ln\frac{rs}{v^2} +\biggl|\ln\frac{r}s\biggl|\biggr)\ln\frac{r}s,\quad
\mathfrak{E}=\biggl(\ln\frac{rs}{v^2} -\biggl|\ln\frac{r}s\biggl|\biggr)\ln\frac{r}s.
\end{gather*}
Then the function $P_{r,s;u,v}(t)$ has the following properties:
\begin{enumerate}
\item
The function $P_{r,s;u,v}(t)$ is increasing on $(0,\infty)$ if and only if either $u>v$, $\mathfrak{A}\ge0$ and $\mathfrak{C}\ge0$ or $u<v$, $\mathfrak{A}\le0$ and $\mathfrak{B}\le0$.
\item
The function $P_{r,s;u,v}(t)$ is decreasing on $(0,\infty)$ if and only if either $u<v$, $\mathfrak{A}\ge0$ and $\mathfrak{B}\ge0$ or $u>v$, $\mathfrak{A}\le0$ and $\mathfrak{C}\le0$.
\item
The function $P_{r,s;u,v}(t)$ is increasing on $(-\infty,0)$ if and only if $u>v$, $\mathfrak{A}\ge0$ and $\mathfrak{E}\ge0$, or $u<v$, $\mathfrak{A}\le0$ and $\mathfrak{D}\le0$.
\item
The function $P_{r,s;u,v}(t)$ is decreasing on $(-\infty,0)$ if and only if either $u>v$, $\mathfrak{A}\le0$ and $\mathfrak{E}\le0$ or $u<v$, $\mathfrak{A}\ge0$ and $\mathfrak{D}\ge0$.
\item
The function $P_{r,s;u,v}(t)$ is increasing on $(-\infty,\infty)$ if and only if either $u>v$, $\mathfrak{C}\ge0$ and $\mathfrak{E}\ge0$ or $u<v$, $\mathfrak{B}\le0$ and $\mathfrak{D}\le0$.
\item
The function $P_{r,s;u,v}(t)$ is decreasing on $(-\infty,\infty)$ if and only if either $u>v$, $\mathfrak{C}\le0$ and $\mathfrak{E}\le0$ or $u<v$, $\mathfrak{B}\ge0$ and $\mathfrak{D}\ge0$.
\item
The function $P_{r,s;u,v}(t)$ on $(-\infty,\infty)$ is logarithmically convex if $\frac{\ln(r/s)}{\ln(u/v)}>1$ or logarithmically concave if $0<\frac{\ln(r/s)}{\ln(u/v)}<1$.
\item
The function $P_{r,s;u,v}(t)$ is $3$-log-convex on $(0,\infty)$ and $3$-log-concave on $(-\infty,0)$ if $\frac{u}v>\frac{r}s>1$ or $\frac{r}s<\frac{u}v<1$; The function $P_{r,s;u,v}(t)$ is $3$-log-concave on $(0,\infty)$ and $3$-log-convex on $(-\infty,0)$ if $\frac{r}s>\frac{u}v>1$ or $\frac{u}v<\frac{r}s<1$.
\end{enumerate}
\end{thm}

\begin{rem}
The monotonicity of the functions $H_{\alpha,\beta;\lambda,\mu}(t)$ and $P_{r,s;u,v}(t)$ can be described by Table~\ref{table1}.
\begin{sidewaystable}
\centering
\label{table1} \caption{Monotoncity of the functions $H_{\alpha,\beta;\lambda,\mu}(t)$ and $P_{r,s;u,v}(t)$}
\begin{tabular}{||c|c||c|c|c|c|c|c||}
  \hline  \hline
  Intervals & Monotonicities & $\mathcal{A}$ or $\mathfrak{A}$ & $\mathcal{B}$ or $\mathfrak{B}$& $\mathcal{C}$ or $\mathfrak{C}$ & $\mathcal{D}$ or $\mathfrak{D}$ & $\mathcal{E}$ or $\mathfrak{E}$ & $\lambda$ and $\mu$ or $u$ and $v$\\ \hline \hline
  $(0,\infty)$ & increasing & $\ge0$ & & $\ge0$ & & & $\lambda>\mu$ or $u>v$ \\  \hline
  $(0,\infty)$ & increasing & $\le0$ & $\le0$ & & & & $\lambda<\mu$ or $u<v$ \\  \hline
  $(0,\infty)$ & decreasing & $\ge0$ & $\ge0$ & & & & $\lambda<\mu$ or $u<v$ \\  \hline
  $(0,\infty)$ & decreasing & $\le0$ & & $\le0$ & & & $\lambda>\mu$ or $u>v$ \\  \hline
  $(-\infty,0)$ & increasing & $\ge0$ & & & & $\ge0$ &$\lambda>\mu$ or $u>v$ \\  \hline
  $(-\infty,0)$ & increasing & $\le0$ & & & $\le0$ & &$\lambda<\mu$ or $u<v$ \\  \hline
  $(-\infty,0)$ & decreasing & $\le0$ & & & & $\le0$ &$\lambda>\mu$ or $u>v$ \\  \hline
  $(-\infty,0)$ & decreasing & $\ge0$ & & & $\ge0$ & &$\lambda<\mu$ or $u<v$ \\  \hline
  $(-\infty,\infty)$ & increasing & & &$\ge0$& &$\ge0$&$\lambda>\mu$ or $u>v$ \\  \hline
  $(-\infty,\infty)$ & increasing &&$\le0$& & $\le0$&&$\lambda<\mu$ or $u<v$ \\  \hline
  $(-\infty,\infty)$ & decreasing && & $\le0$ && $\le0$ &$\lambda>\mu$ or $u>v$ \\ \hline
  $(-\infty,\infty)$ & decreasing && $\ge0$ && $\ge0$ &&$\lambda<\mu$ or $u<v$ \\ \hline\hline
\end{tabular}
\end{sidewaystable}
\normalsize
\end{rem}

\begin{rem}
In~\cite[Remark~2.1]{comp-mon-element-exp.tex} it was remarked that the function $Q_{\alpha,\beta}(t)$ can not be either $4$-log-convex or $4$-log-concave in either $(-\infty,0)$ or $(0,\infty)$, saying nothing of $(-\infty,\infty)$. Therefore, neither $H_{\alpha,\beta;\lambda,\mu}(t)$ nor $P_{r,s;u,v}(t)$ is either $4$-log-convex or $4$-log-concave on either $(-\infty,0)$ or $(0,\infty)$, saying nothing of $(-\infty,\infty)$.
\end{rem}

\section{Proofs of theorems}

\begin{proof}[Proof of Theorem~\ref{exp-gen-T1}]
For $t\ne0$, the function $H_{\alpha,\beta;\lambda,\mu}(t)$ can rewritten as
\begin{equation*}
H_{\alpha,\beta;\lambda,\mu}(t)=\frac{e^{(\alpha-\lambda)t} -e^{(\beta-\lambda)t}}{1-e^{(\mu-\lambda)t}} =\frac{e^{-\frac{\alpha-\lambda}{\mu-\lambda}w} -e^{-\frac{\beta-\lambda}{\mu-\lambda}w}}{1-e^{-w}} =\frac{e^{-Aw}-e^{-Bw}}{1-e^{-w}},
\end{equation*}
where
\begin{equation}
\begin{aligned}
A&=\frac{\alpha-\lambda}{\mu-\lambda},& B&=\frac{\beta-\lambda}{\mu-\lambda}, &w=(\lambda-\mu)t.
\end{aligned}
\end{equation}
Differentiating with respect to $t$ yields
\begin{gather}
H_{\alpha,\beta;\lambda,\mu}'(t)=(\lambda-\mu)Q_{A,B}'(w),\label{h-1-der}\\
\begin{gathered}\label{h-2-der}
[\ln H_{\alpha,\beta;\lambda,\mu}(t)]'' =\biggl[\frac{H_{\alpha,\beta;\lambda,\mu}'(t)} {H_{\alpha,\beta;\lambda,\mu}(t)}\biggr]'
=(\lambda-\mu)\frac{\td}{\td t}\biggl[\frac{Q_{A,B}'(w)}{Q_{A,B}(w)}\biggr]\\ =(\lambda-\mu)^2\biggl[\frac{Q_{A,B}'(w)}{Q_{A,B}(w)}\biggr]'
=(\lambda-\mu)^2[\ln Q_{A,B}(w)]''
\end{gathered}
\end{gather}
and
\begin{equation}\label{h-3-der}
[\ln H_{\alpha,\beta;\lambda,\mu}(t)]''' =(\lambda-\mu)^3\biggl[\frac{Q_{A,B}'(w)}{Q_{A,B}(w)}\biggr]''
=(\lambda-\mu)^3[\ln Q_{A,B}(w)]'''.
\end{equation}
\par
By virtue of~\cite[Theorem~2]{mon-element-exp.tex-rgmia} (see also~\cite[Theorem~3.1]{mon-element-exp-final.tex} and~\cite[Lemma~1]{notes-best.tex}) and the second order derivative~\eqref{h-2-der}, it is easy to deduce that the function $H_{\alpha,\beta;\lambda,\mu}(t)$ is logarithmically convex if $\frac{\beta-\alpha}{\mu-\lambda}>1$ and logarithmically concave if $0<\frac{\beta-\alpha}{\mu-\lambda}<1$ on $(-\infty,\infty)$.
\par
By virtue of~\cite[Theorem~1.1]{comp-mon-element-exp.tex} and the third order derivative~\eqref{h-3-der}, it is not difficult to obtain that
\begin{enumerate}
\item
if $\lambda>\mu$ and $1>\frac{\beta-\alpha}{\mu-\lambda}>0$, then $H_{\alpha,\beta;\lambda,\mu}(t)$ is $3$-log-convex on $(0,\infty)$ and $3$-log-concave on $(-\infty,0)$;
\item
if $\lambda>\mu$ and $\frac{\beta-\alpha} {\mu-\lambda}>1$, then $H_{\alpha,\beta;\lambda,\mu}(t)$ is $3$-log-concave on $(0,\infty)$ and $3$-log-convex on $(-\infty,0)$;
\item
if $\lambda<\mu$ and $1>\frac{\beta-\alpha}{\mu-\lambda}>0$, then $H_{\alpha,\beta;\lambda,\mu}(t)$ is $3$-log-concave on $(0,\infty)$ and $3$-log-convex on $(-\infty,0)$;
\item
if $\lambda<\mu$ and $\frac{\beta-\alpha} {\mu-\lambda}>1$, then $H_{\alpha,\beta;\lambda,\mu}(t)$ is $3$-log-convex on $(0,\infty)$ and $3$-log-concave on $(-\infty,0)$.
\end{enumerate}
\par
Direct computation gives
\begin{align*}
(B-A)(1-A-B)&=\frac{\mathcal{A}}{(\lambda-\mu)^2};\\
(B-A)(|A-B|-A-B)&=
\begin{cases}
\dfrac{\mathcal{B}}{(\lambda-\mu)^2},& \lambda<\mu,\\[1em]
\dfrac{\mathcal{C}}{(\lambda-\mu)^2},& \lambda>\mu;
\end{cases}\\
(B-A)(2-|A-B|-A-B)&=
\begin{cases}
\dfrac{\mathcal{D}}{(\lambda-\mu)^2},& \lambda<\mu,\\[1em]
\dfrac{\mathcal{E}}{(\lambda-\mu)^2},& \lambda>\mu.
\end{cases}
\end{align*}
Consequently, utilization of~\cite[Theorem~2.3]{mon-element-exp-final.tex}, \cite[Corallary~1]{mon-element-exp.tex-rgmia} and the first order derivative~\eqref{h-1-der} yields the following conclusions:
\begin{enumerate}
\item
The function $H_{\alpha,\beta;\lambda,\mu}(t)$ is increasing on $(0,\infty)$ if and only if $\lambda>\mu$, $\mathcal{A}\ge0$ and $\mathcal{C}\ge0$; The function $H_{\alpha,\beta;\lambda,\mu}(t)$ is decreasing on $(0,\infty)$ if and only if $\lambda<\mu$, $\mathcal{A}\ge0$ and $\mathcal{B}\ge0$.
\item
The function $H_{\alpha,\beta;\lambda,\mu}(t)$ is decreasing on $(0,\infty)$ if and only if $\lambda>\mu$, $\mathcal{A}\le0$ and $\mathcal{C}\le0$; The function $H_{\alpha,\beta;\lambda,\mu}(t)$ is increasing on $(0,\infty)$ if and only if $\lambda<\mu$, $\mathcal{A}\le0$ and $\mathcal{B}\le0$.
\item
The function $H_{\alpha,\beta;\lambda,\mu}(t)$ is increasing on $(-\infty,0)$ if and only if $\lambda>\mu$, $\mathcal{A}\ge0$ and $\mathcal{E}\ge0$; The function $H_{\alpha,\beta;\lambda,\mu}(t)$ is decreasing on $(-\infty,0)$ if and only if $\lambda<\mu$, $\mathcal{A}\ge0$ and $\mathcal{D}\ge0$.
\item
The function $H_{\alpha,\beta;\lambda,\mu}(t)$ is decreasing on $(-\infty,0)$ if and only if $\lambda>\mu$, $\mathcal{A}\le0$ and $\mathcal{E}\le0$; The function $H_{\alpha,\beta;\lambda,\mu}(t)$ is increasing on $(-\infty,0)$ if and only if $\lambda<\mu$, $\mathcal{A}\le0$ and $\mathcal{D}\le0$.
\item
The function $H_{\alpha,\beta;\lambda,\mu}(t)$ is increasing on $(-\infty,\infty)$ if and only if $\lambda>\mu$, $\mathcal{C}\ge0$ and $\mathcal{E}\ge0$; The function $H_{\alpha,\beta;\lambda,\mu}(t)$ is decreasing on $(-\infty,\infty)$ if and only if $\lambda<\mu$, $\mathcal{B}\ge0$ and $\mathcal{D}\ge0$.
\item
The function $H_{\alpha,\beta;\lambda,\mu}(t)$ is decreasing on $(-\infty,\infty)$ if and only if $\lambda>\mu$, $\mathcal{C}\le0$ and $\mathcal{E}\le0$. The function $H_{\alpha,\beta;\lambda,\mu}(t)$ is increasing on $(-\infty,\infty)$ if and only if $\lambda<\mu$, $\mathcal{B}\le0$ and $\mathcal{D}\le0$.
\end{enumerate}
The proof of Theorem~\ref{exp-gen-T1} is complete.
\end{proof}

\begin{proof}[Proof of Theorem~\ref{exp-gen-T2}]
This follows directly from the combination of Theorem~\ref{exp-gen-T1} with equations in~\eqref{ratio-relation-1}. Theorem~\ref{exp-gen-T2} is proved.
\end{proof}


\begin{thebibliography}{99}

\bibitem{bullen-handbook}
P. S. Bullen, \textit{Handbook of Means and Their Inequalities}, Mathematics and its Applications (Dordrecht) \textbf{560}, Kluwer Academic Publishers, Dordrecht, 2003.

\bibitem{Best-Constant-exponential.tex}
Ch.-P. Chen and F. Qi, \textit{Best constant in an inequality connected with exponential functions}, Octogon Math. Mag. \textbf{12} (2004), no.~2, 736\nobreakdash--737.

\bibitem{Cheung-Qi-mean-rgmia}
W.-S. Cheung and F. Qi, \textit{Logarithmic convexity of the one-parameter mean values}, RGMIA Res. Rep. Coll. \textbf{7} (2004), no.~2, Art.~15, 331\nobreakdash--342; Available online at \url{http://www.staff.vu.edu.au/rgmia/v7n2.asp}.

\bibitem{Cheung-Qi-mean}
W.-S. Cheung and F. Qi, \textit{Logarithmic convexity of the one-parameter mean values}, Taiwanese J. Math. \textbf{11} (2007), no.~1, 231\nobreakdash--237.

\bibitem{souza}
P. N. de Souza and J.-N. Silva, \textit{Berkeley Problems in Mathematics}, 2nd ed., Problem Books in Mathematics, Springer, New York, 2001.

\bibitem{egp}
N. Elezovi\'c, C. Giordano and J. Pe\v{c}ari\'c, \textit{The best bounds in Gautschi's inequality}, Math. Inequal. Appl. \textbf{3} (2000), 239\nobreakdash--252.

\bibitem{Gauchman-Steffensen-pairs}
H. Gauchman, \textit{Steffensen pairs and associated inequalities}, J. Inequal. Appl. \textbf{5} (2000), no.~1, 53\nobreakdash--61.

\bibitem{best-constant-one-simple.tex}
B.-N. Guo, A.-Q. Liu and F. Qi, \textit{Monotonicity and logarithmic convexity of three functions involving exponential function}, J. Korea Soc. Math. Educ. Ser. B Pure Appl. Math. \textbf{15} (2008), no.~4, 387\nobreakdash--392.

\bibitem{emv-log-convex-simple.tex}
B.-N. Guo and F. Qi, \textit{A simple proof of logarithmic convexity of extended mean values}, Numer. Algorithms (2009), in press; Available online at \url{http://dx.doi.org/10.1007/s11075-008-9259-7}.

\bibitem{ijmest-bernoulli}
B.-N. Guo and F. Qi, \textit{Generalization of Bernoulli polynomials}, Internat. J. Math. Ed. Sci. Tech. \textbf{33} (2002), no.~3, 428\nobreakdash--431.

\bibitem{mon-element-exp-final.tex}
B.-N. Guo and F. Qi, \textit{Properties and applications of a function involving exponential functions}, Commun. Pure Appl. Anal. \textbf{8} (2009), no.~4, 1231\nobreakdash--1249.

\bibitem{remiander-Sen-Lin-Guo.tex}
S. Guo and F. Qi, \textit{A class of completely monotonic functions related to the remainder of Binet's formula with applications}, Tamsui Oxf. J. Math. Sci. \textbf{25} (2009), no.~1, in press.

\bibitem{3rded}
J.-Ch. Kuang, \textit{Ch\'angy\`ong B\`ud\v{e}ngsh\`\i} (\textit{Applied Inequalities}), 3rd ed., Shandong Science and Technology Press, Ji'nan City, Shandong Province, China, 2004. (Chinese)

\bibitem{best-constant-one.tex}
A.-Q. Liu, G.-F. Li, B.-N. Guo and F. Qi, \textit{Monotonicity and logarithmic concavity of two functions involving exponential function}, Internat. J. Math. Ed. Sci. Tech. \textbf{39} (2008), no.~5, 686\nobreakdash--691.

\bibitem{bernoulli-luo-guo-qi-rgmia}
Q.-M. Luo, B.-N. Guo and F. Qi, \textit{Generalizations of Bernoulli's numbers and polynomials}, RGMIA Res. Rep. Coll. \textbf{5} (2002), no.~2, Art.~12, 353\nobreakdash--359; Available online at \url{http://www.staff.vu.edu.au/rgmia/v5n2.asp}.

\bibitem{bernoulli-luo-guo-qi-debnath-IJMMS}
Q.-M. Luo, B.-N. Guo, F. Qi, and L. Debnath, \textit{Generalizations of Bernoulli numbers and polynomials}, Internat. J. Math. Math. Sci. \textbf{2003} (2003), no.~59, 3769\nobreakdash--3776.

\bibitem{euler-bernoulli-luo-qi-adv}
Q.-M. Luo and F. Qi, \textit{Relationships between generalized Bernoulli numbers and polynomials and generalized Euler numbers and polynomials}, Adv. Stud. Contemp. Math. (Kyungshang) \textbf{7} (2003), no.~1, 11\nobreakdash--18.

\bibitem{sandor-gamma-3.tex}
F. Qi, \textit{A class of logarithmically completely monotonic functions and the best bounds in the first Kershaw's double inequality}, J. Comput. Appl. Math. \textbf{206} (2007), no.~2, 1007\nobreakdash--1014.

\bibitem{best-constant-rgmia}
F. Qi, \textit{A monotonicity result of a function involving the exponential function and an application}, RGMIA Res. Rep. Coll. \textbf{7} (2004), no.~3, Art.~16, 507\nobreakdash--509; Available online at \url{http://www.staff.vu.edu.au/rgmia/v7n3.asp}.

\bibitem{schext}
F. Qi, \textit{A note on Schur-convexity of extended mean values}, Rocky Mountain J. Math. \textbf{35} (2005), no.~5, 1787\nobreakdash--1793.

\bibitem{Wendel-Gautschi-type-ineq.tex}
F. Qi, \textit{Bounds for the ratio of two gamma functions---From Wendel's and related inequalities to logarithmically completely monotonic functions}, Available online at \url{http://arxiv.org/abs/0904.1048}.

\bibitem{Wendel2Elezovic.tex}
F. Qi, \textit{Bounds for the ratio of two gamma functions---From Wendel's limit to Elezovi\'c-Giordano-Pe\v{c}ari\'c's theorem}, Available online at \url{http://arxiv.org/abs/0902.2514}.

\bibitem{pams-62}
F. Qi, \textit{Logarithmic convexity of extended mean values}, Proc. Amer. Math. Soc. \textbf{130} (2002), no.~6, 1787\nobreakdash--1796.

\bibitem{pams-62-rgmia}
F. Qi, \textit{Logarithmic convexities of the extended mean values}, RGMIA Res. Rep. Coll. \textbf{2} (1999), no.~5, Art.~5, 643\nobreakdash--652; Available online at \url{http://www.staff.vu.edu.au/rgmia/v2n5.asp}.

\bibitem{mon-element-exp.tex-rgmia}
F. Qi, \textit{Monotonicity and logarithmic convexity for a class of elementary functions involving the exponential function}, RGMIA Res. Rep. Coll. \textbf{9} (2006), no.~3, Art.~3; Available online at \url{http://www.staff.vu.edu.au/rgmia/v9n3.asp}.

\bibitem{schext-rgmia}
F. Qi, \textit{Schur-convexity of the extended mean values}, RGMIA Res. Rep. Coll. \textbf{4} (2001), no.~4, Art.~4, 529\nobreakdash--533; Available online at \url{http://www.staff.vu.edu.au/rgmia/v4n4.asp}.

\bibitem{cubo}
F. Qi, \textit{The extended mean values: Definition, properties, monotonicities, comparison, convexities, generalizations, and applications}, Cubo Mat. Educ. \textbf{5} (2003), no.~3, 63\nobreakdash--90.

\bibitem{comp-mon-element-exp.tex}
F. Qi, \textit{Three-log-convexity for a class of elementary functions involving exponential function}, J. Math. Anal. Approx. Theory \textbf{1} (2006), no.~2, 100\nobreakdash--103.

\bibitem{Cheung-Qi-Rev.tex}
F. Qi, P. Cerone, S. S. Dragomir and H. M. Srivastava, \textit{Alternative proofs for monotonic and logarithmically convex properties of one-parameter mean values}, Appl. Math. Comput. \textbf{208} (2009), no.~1, 129\nobreakdash--133; Available online at \url{http://dx.doi.org/10.1016/j.amc.2008.11.023}.

\bibitem{steffensen-qi-cheng-rgmia}
F. Qi and J.-X. Cheng, \textit{New Steffensen pairs}, RGMIA Res. Rep. Coll. \textbf{3} (2000), no.~3, Art.~11, 431\nobreakdash--436; Available online at \url{http://www.staff.vu.edu.au/rgmia/v3n3.asp}.

\bibitem{steffensen-pair-Anal}
F. Qi and J.-X. Cheng, \textit{Some new Steffensen pairs}, Anal. Math. \textbf{29} (2003), no.~3, 219\nobreakdash--226.

\bibitem{qcw}
F. Qi, J.-X. Cheng and G. Wang, \textit{New Steffensen pairs}, Inequality Theory and Applications, Volume \textbf{1}, 273\nobreakdash--279, Nova Science Publishers, Huntington, NY, 2001.

\bibitem{notes-best-new-proof.tex}
F. Qi and B.-N. Guo, \textit{An alternative proof of Elezovi\'c-Giordano-Pe\v{c}ari\'c's theorem}, Available online at \url{http://arxiv.org/abs/0903.1174}.

\bibitem{bernoulli-qi-guo-rgmia}
F. Qi and B.-N. Guo, \textit{Generalisation of Bernoulli polynomials}, RGMIA Res. Rep. Coll. \textbf{4} (2001), no.~4, Art.~10, 691\nobreakdash--695; Available online at \url{http://www.staff.vu.edu.au/rgmia/v4n4.asp}.

\bibitem{onsp}
F. Qi and B.-N. Guo, \textit{On Steffensen pairs}, J. Math. Anal. Appl. \textbf{271} (2002), no.~2, 534\nobreakdash--541.

\bibitem{onsp-rgmia}
F. Qi and B.-N. Guo, \textit{On Steffensen pairs}, RGMIA Res. Rep. Coll. \textbf{3} (2000), no.~3, Art.~10, 425\nobreakdash--430; Available online at \url{http://www.staff.vu.edu.au/rgmia/v3n3.asp}.

\bibitem{exp-funct-further.tex}
F. Qi and B.-N. Guo, \textit{The function $(b^x-a^x)/x$: Logarithmic convexity}, RGMIA Res. Rep. Coll. \textbf{11} (2008), no.~1, Art.~5; Available online at \url{http://www.staff.vu.edu.au/rgmia/v11n1.asp}.

\bibitem{exp-funct-appl-means.tex}
F. Qi and B.-N. Guo, \textit{The function $(b^x-a^x)/x$: Logarithmic convexity and applications to extended mean values}, Available online at \url{http://arxiv.org/abs/0903.1203}.

\bibitem{sandor-gamma-3-note.tex-aam}
F. Qi and B.-N. Guo, \textit{Wendel's and Gautschi's inequalities: Refinements, extensions, and a class of logarithmically completely monotonic functions}, Appl. Math. Comput. \textbf{205} (2008), no.~1, 281\nobreakdash--290; Available online at \url{http://dx.doi.org/10.1016/j.amc.2008.07.005}.

\bibitem{sandor-gamma-3-note.tex}
F. Qi and B.-N. Guo, \textit{Wendel-Gautschi-Kershaw's inequalities and sufficient and necessary conditions that a class of functions involving ratio of gamma functions are logarithmically completely monotonic}, RGMIA Res. Rep. Coll. \textbf{10} (2007), no.~1, Art.~2; Available online at \url{http://www.staff.vu.edu.au/rgmia/v10n1.asp}.

\bibitem{notes-best.tex}
F. Qi, B.-N. Guo and Ch.-P. Chen, \textit{The best bounds in Gautschi-Kershaw inequalities}, Math. Inequal. Appl. \textbf{9} (2006), no.~3, 427\nobreakdash--436.

\bibitem{note-on-neuman.tex}
F. Qi, S. Guo and B.-N. Guo, \textit{A class of $k$-log-convex functions and their applications to some special functions}, RGMIA Res. Rep. Coll. \textbf{10} (2007), no.~1, Art.~21; Available online at \url{http://www.staff.vu.edu.au/rgmia/v10n1.asp}.

\bibitem{note-on-neuman-ITSF-simplified.tex}
F. Qi, S. Guo, B.-N. Guo and Sh.-X. Chen, \textit{A class of $k$-log-convex functions and their applications to some special functions}, Integral Transforms Spec. Funct. \textbf{19} (2008), no.~3, 195\nobreakdash--200.

\bibitem{ql}
F. Qi and Q.-M. Luo, \textit{A simple proof of monotonicity for extended mean values}, J. Math. Anal. Appl. \textbf{224} (1998), no.~2, 356\nobreakdash--359.

\bibitem{psi-reminders.tex}
F. Qi, D.-W. Niu and B.-N. Guo, \textit{Monotonic properties of differences for remainders of psi function}, Internat. J. Pure Appl. Math. Sci. \textbf{4} (2007), no.~1, ???--???.

\bibitem{psi-reminders.tex-rgmia}
F. Qi, D.-W. Niu and B.-N. Guo, \textit{Monotonic properties of differences for remainders of psi function}, RGMIA Res. Rep. Coll. \textbf{8} (2005), no.~4, Art.~16, 683\nobreakdash--690; Available online at \url{http://www.staff.vu.edu.au/rgmia/v8n4.asp}.

\bibitem{jmaa-ii-97}
F. Qi and S.-L. Xu, \textit{Refinements and extensions of an inequality, I\!I}, J. Math. Anal. Appl. \textbf{211} (1997), no.~2, 616\nobreakdash--620.

\bibitem{(b^x-a^x)/x}
F. Qi and S.-L. Xu, \textit{The function $(b^x-a^x)/x$: Inequalities and properties}, Proc. Amer. Math. Soc. \textbf{126} (1998), no.~11, 3355\nobreakdash--3359.

\bibitem{qx3}
F. Qi, S.-L. Xu and L. Debnath, \textit{A new proof of monotonicity for extended mean values}, Internat. J. Math. Math. Sci. \textbf{22} (1999), no.~2, 415\nobreakdash--420.

\bibitem{wuzh}
Y.-D. Wu and Zh.-H. Zhang, \textit{The best constant for an inequality}, Octogon Math. Mag. \textbf{12} (2004), no.~1, 139\nobreakdash--141.

\bibitem{wuzh-rgmia}
Y.-D. Wu and Zh.-H. Zhang, \textit{The best constant for an inequality}, RGMIA Res. Rep. Coll. \textbf{7} (2004), no.~1, Art.~19; Available online at \url{http://www.staff.vu.edu.au/rgmia/v7n1.asp}.

\bibitem{best-constant-one-simple-real.tex}
Sh.-Q. Zhang, B.-N. Guo and F. Qi, \textit{A concise proof for properties of three functions involving the exponential function}, Appl. Math. E-Notes \textbf{9} (2009), in press.

\end{thebibliography}
\end{document}